\documentclass[11pt]{amsart}
\usepackage{amscd}
\usepackage{amsmath}
\usepackage{graphicx}
\usepackage{amsfonts}
\usepackage{amssymb}
\newtheorem{theorem}{Theorem}[section]

\newtheorem{corollary}[theorem]{Corollary}
\newtheorem{lemma}[theorem]{Lemma}
\newtheorem{example}[theorem]{Example}
\newtheorem{remark}[theorem]{Remark}

\def\IR{{\mathbb R}}
\def\IC{{\mathbb C}}

\def\cA{{\mathcal A}}
\def\cB{{\mathcal B}}

\def\cF{{\mathcal F}}
\def\bT{{\mathbb T}}
\def\cS{{\mathcal S}}
\def\qed{\hfill $\Box$\medskip}
\def\BH{{\cB(H)}}

\begin{document}
\title[unitary similarity invariant function preservers]{Unitary
similarity invariant function preservers of skew products of operators}
\author{Jianlian Cui}
\address[Jianlian Cui]{
Department of Mathematical Science, Tsinghua University, Beijing
100084, P.R. China.} \email{\rm jcui@math.tsinghua.edu.cn}
\author{Chi-Kwong Li}
\address[Chi-Kwong Li] {Department of Mathematics, The College of William
\& Mary, Williamsburg, VA 13185, USA. {\rm Li is an honorary
professor of the University of Hong Kong,
and the Shanghai University.}}
\email{\rm ckli@math.wm.edu}
\author{Nung-Sing Sze}
\address[Nung-Sing Sze]
{Department of Applied Mathematics, The Hong Kong Polytechnic University, Hong Kong.}
\email{\rm raymond.sze@polyu.edu.hk}



\thanks{Research of the first author
was supported by National Natural Science Foundation of
China (No.11271217). Research of the second author was partially
supported by the USA NSF DMS 1331021, the Simons Foundation Grant 351047,
and NNSF of China Grant 11571220.
Research of the third author was supported by a HK RGC grant PolyU
502512 and a PolyU central research grant G-YBCR}

\maketitle

\begin{abstract}

Let ${\mathcal B}(H)$ denote the Banach algebra of all bounded
linear operators on a complex Hilbert space $H$ with
$\dim H\geq 3$, and let $\mathcal A$ and $\mathcal B$ be subsets of
${\mathcal B}(H)$ which contain all rank one operators.
Suppose $F(\cdot )$ is
a unitary invariant norm, the pseudo spectra, the pseudo spectral radius,
the $C$-numerical range, or the $C$-numerical radius for some finite rank operator $C$.
The structure is determined for surjective maps $\Phi :{\mathcal A}\rightarrow \mathcal B$
satisfying $F(A^*B)=F(\Phi (A)^*\Phi (B))$ for all $A, B \in {\mathcal A}$.
To establish the proofs, some general results are obtained for functions
$F:{\mathcal F}_1(H) \cup \{0\} \rightarrow [0, +\infty)$, where ${\mathcal F}_1(H)$
is the set of rank one operators in ${\mathcal B}(H)$,
satisfying 
(a) $F(\mu UAU^*)=F(A)$ for a complex unit $\mu$, $A\in {\mathcal F}_1(H)$ and unitary $U \in {\mathcal B}(H)$,
(b) for any rank one operator $X\in {\mathcal F}_1(H)$ the map $t\mapsto F(tX)$ on $[0, \infty)$ is strictly increasing, and
(c) the set $\{F(X): X \in {\mathcal F}_1(H) \hbox{ and } \|X\| = 1\}$ attains its maximum and minimum.
\end{abstract}

\section{Introduction}

There has been considerable interest in studying maps $\Phi$ on matrices or operators satisfying
$F(\Phi(A)\diamond \Phi(B)) = F(A\diamond B)$ for different kind of functions $F(\cdot)$ and
different kinds of product $A \diamond B$. We say that such a map
preserves the function $F$ of the product $A \diamond B$ on operators.
For example, researchers have considered such problems for functions $F$ including the
spectrum, the spectral radius, a unitary invariant norm, a unitary similarity invariant norm,
a generalized numerical range, and numerical radius for the product $A\diamond B$ such as the
usual product $AB$,  the Lie product $AB-BA$, the Jordan product $AB+BA$, the Jordan triple product,
and the Schur (entrywise) product on matrices;
for example, see \cite{CLS,CLS2,CL,CL1,CH1,CH2,GL,HHZ1,HLW1,HLW2,LP,LPR,LPS,LS} and their references.

In this paper, we consider maps preserving
$F$ for the  skew product $A \diamond B = A^*B$ of operators, where $F$
is a unitary invariant norm, a unitary similarity invariant norm, the pseudo spectrum,
the pseudo spectral radius, and the $C$-numerical radius.
We obtain a general result for some general function $F$ so that
one can use the general result to treat the special cases.

Denote by $\mathbb N$, $\mathbb C$ and $\mathbb T$
the set of natural numbers, the complex field and the unit circle of $\mathbb C$, respectively.
Let ${\mathcal B}(H)$ be the Banach algebra of
all bounded linear operators on a complex Hilbert space $H$, and let
${\mathcal F}(H)$ and ${\mathcal F}_1(H)$ be the sets of
all finite rank linear operators 
and 
all rank one linear operators 
in ${\mathcal B}(H)$, respectively.
For any
$x$, $f\in H$, the notation $x\otimes f$ denotes a rank one operator
on $H$ defined by $z\mapsto \langle z,f\rangle x$ for every $z\in
H$; and every operator of rank one in ${\mathcal B}(H)$ can be
written in this form. 
Fix an arbitrary orthogonal basis
$\{e_i\}_{i\in \Gamma }$ of $H$. For $x\in H$, write
$x=\sum_{i\in \Gamma }\xi _ie_i$, and define
the conjugate operator $J: H \to H$ by $J  x= \bar{x}
=\sum _{i\in \Gamma}\bar{\xi _i}e_i$.
Finally, the notation $\overline A$ denotes the bounded linear operator $JAJ$ in $\cB(H)$.
Notice that $\langle  \overline A e_i, e_j \rangle = \overline{ \langle A e_i, e_j \rangle}$
for all $i,j \in \Gamma$.

\section{Unitary similarity invariant functions preservers}

In this section, we consider maps preserving unitary similarity invariant functions
of skew product of operators.  In particular, we consider a function
$F:{\cF_1(H)}\cup \{0\} \rightarrow [0, +\infty)$ satisfying some of the following
properties:
\begin{enumerate}
\item[\rm (F1)] $F(\mu UXU^*)=F(X)$ for any complex unit $\mu$, $X\in \cF_1(H)$
and unitary $U \in \BH$.


\item[\rm  (F2)] For any  $X\in \cF_1(H)$, the map $t\mapsto F(tX)$ on $[0, \infty)$ is strictly increasing.

\item[\rm (F3)] The set $\{F(X): X \in \cF_1(H) \hbox{ and } \|X\| = 1\}$ attains its maximum and minimum.
\end{enumerate}

\medskip
Notice that (F2) implies $F(X) > 0$ for all $X \in \cF_1(H)$.
We will prove the following.

\begin{theorem}\label{main}
Let $H$ be a complex Hilbert space with
$\dim H\geq 3$, $\mathcal A$ and $\mathcal B$ be subsets of
${\mathcal B}(H)$ which contain $\cF_1(H)$.
Suppose $\Phi :{\cA}\rightarrow {\cB}$ is a surjective map such that
\begin{equation}\label{zero}
A^*B = 0 \quad\Longleftrightarrow \quad \Phi(A)^* \Phi(B) = 0 \quad \hbox{for all}\quad A,B \in \cA.
\end{equation}
Then $\Phi$ preserves rank one operators in both directions, and $\Phi(0) = 0$, if $0 \in \cA$.
Moreover,  there exist a unitary operator $U$
in ${\mathcal B}(H)$ and a map $h: H \times H \to H$ 
such that
$$\Phi(x \otimes f) = U x \otimes h(x,f)
\quad\hbox{for all}\quad x,f\in H,$$
or
$$\Phi(x \otimes f) = U J x \otimes h(x,f)
\quad\hbox{for all}\quad x,f\in H.$$
\end{theorem}

\begin{theorem} \label{main-2} Let $H, \cA, \cB$ and $\Phi$ satisfy the 
hypotheses of Theorem {\ref{main}}.
Assume that $F:{\cF_1(H)}\cup \{0\} \rightarrow [0, +\infty)$ is a function such that 
\begin{equation}\label{FAB}
F(\Phi(A)^*\Phi(B) ) = F(A^*B) \quad\hbox{whenever $A$ or $B$ has rank one.}
\end{equation}
If {\rm (F1)} and {\rm (F2)} hold, then the map $h$ in Theorem \ref{main} satisfies $\|h(x,f)\| = \|f\|$ for all $x,f \in H$.
If, in addition, {\rm (F3)} is satisfied, then there exist a unitary $U$ and  a partial isometry $V_A$
on ${\mathcal B}(H)$,
where $V_A$ depends on $A$ and $V_A^*V_A$ is the right support projection of $A$,
such that either
\begin{equation}\label{form-1}
\Phi(A) = U AV_A^* \quad\hbox{for every}\quad A \in {\mathcal A},
\end{equation}
or 
\begin{equation}\label{form-2}
\Phi(A) = U JAJ V_A^* \quad\hbox{for every}\quad A \in {\mathcal A}.
\end{equation}
\end{theorem}

\begin{remark} \label{main-3} \rm
Note that in Theorem \ref{main-2}, if  
$F(X) = F(JXJ)$ for all $X\in \cF_1(H)$, we can assume the $\Phi$ has any one of the 
form (\ref{form-1}) or (\ref{form-2}).
\end{remark}

For $x\in H$, let $L_x=\{x\otimes f\mid f\in H\}$.
To prove Theorem \ref{main}, we need an auxiliary result.
Similar to \cite[Theorem 2.1]{HHZ1}, the following lemma can be obtained.
See also \cite[Lemma 2.2]{CH2}.

\begin{lemma} \label{2.2}
Let $H$, $\mathcal A$ and $\mathcal B$ be just as assumptions in Theorem \ref{main}.
Suppose $\Phi :{\mathcal A}\rightarrow {\mathcal B}$ is a surjective map
satisfying (\ref{zero}) in Theorem \ref{main}.
Then $\Phi$ preserves rank one operators in both directions and there exists a unitary
operator $U$ on $H$ such that
$\Phi (L_x)=L_{Ux}$  for every $x\in H$ or
$\Phi (L_x)=L_{UJx}$  for every $x\in H$.
\end{lemma}

The referee suggested that it might be possible to relax the assumption (\ref{zero}) in Lemma \ref{2.2} (and also Theorem \ref{main}) to 
``$\Phi (A)^*\Phi (B)=0$ if and only if $A^*B=0$ for all $A, B\in {\mathcal F}_1(H)$''.
Unfortunately, it is impossible in view of the following example.

\begin{example} \label{example}
Let $H$ be a separable Hilbert space with the orthonormal basis $\{e_k:k \in {\mathbb N} \}$.
Suppose $S$ is the bounded linear operator on $H$ defined by $S e_k = e_{k+1}$ for all $k \in {\mathbb N}$.
Let ${\mathcal R} =  \{x \otimes f + x_1 e_1 \otimes Sf: x,f \in H, x_1 \in \IC\}$
and define the  map 
$\Phi: {\mathcal F}_1(H) \cup {\mathcal R} \to {\mathcal F}_1(H)$  by
$$\Phi(x\otimes f) = Sx \otimes f \quad\hbox{and}\quad
\Phi(x \otimes f + x_1 e_1 \otimes Sf) = (x_1 e_1 + Sx) \otimes f
\quad\hbox{for all } x,f\in H,\, x_1 \in \IC.$$ Then $\Phi$ is surjective and 
$$
A^*B = 0 \quad\Longleftrightarrow \quad \Phi(A)^* \Phi(B) = 0 \quad \hbox{whenever $A$ 
or $B$ has rank one}.
$$
\end{example}

Clearly, $\Phi$ does not satisfy the conclusion stated in Lemma \ref{2.2}.
Therefore, the condition (\ref{zero}) of Theorem \ref{main} cannot be relaxed unless one can find a new proof without using Lemma \ref{2.2}
by the weaker assumption.
Now we are in a position to prove Theorems \ref{main} and \ref{main-2}.

\medskip
{\it Proof of Theorem \ref{main}.}
Let $\Phi:{\mathcal A}\rightarrow {\mathcal B}$ be a surjective map satisfying (\ref{zero}).
It follows from Lemma \ref{2.2} that $\Phi$ preserves rank one operators in both directions
and there exists a unitary
operator $U$ on $H$ such that
$\Phi (L_x)=L_{Ux}$ for every $x\in H$
or
$\Phi (L_x)=L_{UJx}$ for every $x\in H$.
Notice that $A^*B = 0$ for all $B \in \cF_1(H)$ if and only if $A = 0$.
Therefore, $\Phi(0) = 0$ if $0\in \cA$.
Finally, from the definition of $L_x$, there exists 
a map $h: H \times H \to H$ 
such that
$\Phi(x \otimes f) = U x \otimes h(x,f)$
for all $x,f\in H$,
or
$\Phi(x \otimes f) = U J x \otimes h(x,f)$
for all $x,f\in H$.
\qed

\medskip
{\it Proof of Theorem \ref{main-2}.}
Suppose $\Phi$ satisfies condition (\ref{zero}) in Theorem \ref{main}. Then $\Phi$ has one of the two forms stated in Theorem \ref{main}.
By replacing $\Phi$ with the map $A \mapsto \Phi(\overline A) = \Phi(JAJ)$ for the latter case,
we may always assume the former case holds.
That is, for any nonzero
$x, f\in H$, there is $h(x, f)\in H$ such that $\Phi (x\otimes f)=Ux\otimes h(x, f)$.
Suppose now a function $F: \cF_1(H) \cup \{0\} \to [0,\infty)$ has the properties (F1) and (F2) and $\Phi$ also satisfies (\ref{FAB}).
For any $x, f\in H$, the equality
$F(\Phi (x\otimes f)^*\Phi(x\otimes f))=F((x\otimes f)^*x\otimes f)$ ensures that
$$F(\langle x,x \rangle h(x,f)\otimes h(x,f))=F(\langle x,x \rangle  f\otimes f).$$
It follows from conditions (F1) and (F2)  that
$\|h(x, f)\|=\|f\|$ for all $x, f\in H$. Thus, the first result follows.
Furthermore, for any nonzero $\alpha \in \IC$,
$$Ux \otimes h(x,f) = \Phi(x\otimes f) = \Phi(\alpha^{-1} x \otimes \bar{\alpha} f)
= U(\alpha^{-1} x) \otimes h(\alpha^{-1} x, \bar{\alpha }f)
= Ux  \otimes\bar{  \alpha}^{-1} h(\alpha^{-1} x, \bar{\alpha} f),$$
and hence $h(\alpha^{-1} x , \bar{\alpha} f) = \bar{\alpha} h(x,f)$ for all nonzero $\alpha \in \IC$.
For any $x, f\in H$ and $A \in \cA$, we have
\begin{equation}\label{FAU}
F(\Phi (A)^*Ux\otimes h(x,f)) =F(\Phi (A)^*\Phi (x\otimes f))=F(A^*x\otimes f).
\end{equation}
Now we further assume that $F$ also satisfies condition (F3). Together with condition (F1),
this is equivalent to say that for any fixed nonzero $x \in H$,
the set
\begin{equation}\label{ext}
\left\{F(x\otimes f) : f\in H,\,  \|f\| = \|x\|^{-1} \right\} \hbox{ always attains its maximum and minimum points.}
\end{equation}
Clearly, there exists a unitary operator $W$ on $H$ such that
$W\frac{\Phi (A)^*Ux}{\|\Phi (A)^*Ux\|}=\frac{A^*x}{\|A^*x\|}$. Then
\begin{equation}\label{F}
F\left ( \frac{\|\Phi(A)^*U x\|}{\|A^*x\|} A^*x \otimes W h(x,f) \right)
= F\left(\frac{\|\Phi(A)^*U x\|}{\|A^*x\|} W^* A^*x \otimes h(x,f) \right)
= F\left(A^*x \otimes f\right).
\end{equation}
Suppose $\|\Phi(A)^*U x\| > \|A^*x\|$.
Then (\ref{F}) and condition (F3) imply that for any $f\in H$,
$$F\left ( A^*x \otimes W h(x,f) \right) <
F\left ( \frac{\|\Phi(A)^*U x\|}{\|A^*x\|} A^*x \otimes W h(x,f) \right)
= F\left(A^*x \otimes f\right).
$$
Then the  set $\{F(A^*x \otimes f) : f \in H, \|f\| = \|A^*x\|^{-1}\}$ cannot attain its minimum
because $\|Wh(x,f) \| = \|h(x,f)\| = \|f\|$.
This contradicts condition (\ref{ext}).
Thus, $\|\Phi(A)^*Ux\| \le \|A^*x\|$.
Now we assume that $\|\Phi(A)^*Ux\| < \|A^*x\|$.
Again by (\ref{F}) and condition (F3), for any $f \in H$,
$$F\left ( \frac{\|\Phi(A)^*U x\|}{\|A^*x\|} A^*x \otimes W h(x,f) \right)
= F\left(A^*x \otimes f\right)
> F\left(\frac{\|\Phi(A)^*U x\|}{\|A^*x\|} A^*x \otimes f\right).$$
With the fact that $\|Wh(x,f) \| = \|f\|$, one can conclude that
the set $$\left\{ F \left( \frac{\|\Phi(A)^*U x\|}{\|A^*x\|}A^*x \otimes f \right) : f \in H,\, \|f\|
= \|\Phi(A)^*Ux\|^{-1} \right\}$$
cannot attains its maximum,
which again contradicts condition (\ref{ext}).
Thus, one concludes that
\begin{eqnarray}\label{norm}
\|\Phi (A)^*Ux\|=\|A^*x\|\quad \hbox{for every }x\in H.
\end{eqnarray}
This implies that there exists a partial isometry $V_A$
on $\cB(H)$ such that $V_A^* V_A$ is the right support projection of $A$ and
$$V_AA^*x=\Phi (A)^*Ux \quad \hbox{for all } x\in H.$$
So $\Phi(A)=UAV_A^*$.
\qed

\noindent
\begin{remark} \label{remark1} \rm
The same result can also be obtained if one replaces conditions (F2)  and (F3)
in Theorem \ref{main} by the following properties.
\begin{enumerate}
\item[\rm  (F2')] $F(X) > 0$ for all $X \in \cF_1(H)$ and there exists a strictly increasing function $g$ on $[0,\infty)$
such that $F(tX) = g(t) F(X)$ for any $t \in [0,\infty)$ and
any $X \in \cF_1(H)$.

\item[\rm (F3')] The set $\{F(X): X \in \cF_1(H) \hbox{ and } \|X\| = 1\}$  is bounded.
\end{enumerate}
\end{remark}

\noindent
Notice that  condition (F2') is stronger than (F2) while  condition (F3') is weaker than (F3).

\medskip
\it Proof of Remark \ref{remark1}. \rm Suppose  conditions (F2) and (F3)
are replaced by conditions (F2') and (F3') in Theorem \ref{main-2}.
Following the proof of Theorem \ref{main-2} up to equation (\ref{F}),
we conclude that for any $x,f\in H$ and $A \in \cA$,
$$F\left(A^*x \otimes f\right)
= F\left ( s A^*x \otimes W h(x,f) \right)
= g(s) F\left(A^* x \otimes W h(x,f) \right),$$
where $s = \frac{\|\Phi(A)^*U x\|}{\|A^*x\|}$.
It follows from conditions (F1) and (F3') that the set
$\{F(A^* x \otimes f): f\in H,\, \|f\| = \|A^*x\|^{-1}\}$ is always bounded.
Let $R$ and $S$ be the infimum and supremum of this set. Then
$\|Wh(x,f)\| = \|h(x,f)\| = \|f\|$ implies that
$$g(s) R \le F(A^*x \otimes f) = g(s) F\left(A^* x \otimes W h(x,f) \right) \le g(s) S
\quad\Longrightarrow\quad
g(s) R \le R \le S \le g(s)  S.$$
It follows that $g(s) = 1$. Since $g$ is strictly increasing 
and $F(1\!\cdot\! X) = g(1) F(X)$,
$g(s) = 1$ if and only if $s = 1$. Thus, $\|\Phi(A)^*Ux\| = \|A^*x\|$ for all $x \in H$
and hence equation (\ref{norm}) holds.
Then the remaining steps follow from the same argument as in the proof of Theorem \ref{main-2}.
\qed

\section{unitary invariant function preservers}

\medskip
If the function under consideration is unitary invariant, i.e.,
\begin{itemize}
\medskip
\item[\rm (F1'')] $F(UAV) = F(A)$ for any $A \in \cF_1(H)$ and unitary $U, V \in \BH$,
\end{itemize}
\medskip\noindent
then  the set $\{F(X): X \in \cF_1(H) \hbox{ and } \|X\| = 1\}$ is a singleton.
Hence, condition (F3) always holds in this case. Also, condition (F2) reduces to
\begin{enumerate}
\item[\rm  (F2")] For any rank one projection $X$, the map $t\mapsto F(tX)$ on
$[0, \infty)$ is strictly increasing.
\end{enumerate}
In fact, take $f=A^*x$ in equation (\ref{FAU}). Since $\|h(x, f)\|=\|f\|$ for any $x, f\in H$,
there exists a unitary operator $W_1$ on $H$ such that $W_1h(x,f)=f$. Also there exists a
unitary operator $W_2$ on $H$ such that
$W_2\frac{\Phi (A)^*Ux}{\|\Phi (A)^*Ux\|}=\frac{A^*x}{\|A^*x\|}$. Thus, equation (\ref{FAU}),
together with condition (F1''), entails that
$$
F\left (\|\Phi (A)^*Ux\|\|A^*x\|\frac{A^*x}{\|A^*x\|}\otimes \frac{A^*x}{\|A^*x\|} \right) =F\left( \|A^*x\|^2\frac{A^*x}{\|A^*x\|}\otimes \frac{A^*x}{\|A^*x\|} \right),$$
which, together with condition (F2"), implies that
$$\|\Phi (A)^*Ux\|=\|A^*x\|\quad \mbox{ for every }x\in H.$$Now the result
follows from the corresponding proof of Theorem \ref{main-2}.
Thus, we have the following theorem.

\begin{theorem}\label{uif}
Let $H$ be a complex Hilbert space with
$\dim H\geq 3$, $\mathcal A$ and $\mathcal B$ be subsets of
${\mathcal B}(H)$ which contain $\cF_1(H)$.
Suppose $\Phi :{\cA}\rightarrow {\cB}$ is a surjective map satifiying (\ref{zero}) in Theorem \ref{main}.
Assume that $F:{\cF_1(H)}\cup \{0\} \rightarrow [0, +\infty)$ is a function such that 
$$F(\Phi(A)^*\Phi(B) ) = F(A^*B),$$
whenever $A$ or $B$, and thus, $\Phi(A)$ or $\Phi(B)$, has rank one.
If {\rm (F1'')} and {\rm (F2'')} hold, then there exist a unitary $U$ and a partial isometry $V_A$
on ${\mathcal B}(H)$,
where $V_A$ depends on $A$ and $V_A^*V_A$ is the right support projection of $A$,
such that either
\begin{equation}\label{form-a}
\Phi(A) = U AV_A^* \quad\hbox{for every}\quad A \in {\mathcal A},
\end{equation}
or
\begin{equation}\label{form-b}
\Phi(A) = U JAJ V_A^* \quad\hbox{for every}\quad A \in {\mathcal A}.
\end{equation}
\end{theorem}

Similar to Remark \ref{main-3},
in the last assertion of Theorem \ref{uif}
if $F(X) = F(JXJ)$ for all $X\in \cF_1(H)$, then
one may assume that $\Phi$ has any one of the 
form (\ref{form-a}) or (\ref{form-b}).

Recall that a norm $|\!|\!| \cdot |\!|\!|$ on $\BH$ is unitarily invariant if
$|\!|\!| UAV |\!|\!| = |\!|\!|A |\!|\!|$ for all
$A \in \BH$ and unitary $U, V \in \BH$.

\begin{corollary}\label{uin}
Let $H$, $\mathcal A$ and $\mathcal B$ just as assumptions in Theorem \ref{uif}.
Assume that $|\!|\!|\cdot
|\!|\!|$ is a unitary invariant norm. Suppose a surjective map
$\Phi:{\mathcal A}\rightarrow {\mathcal B}$ satisfies
$|\!|\!|\Phi(A)^*\Phi (B) |\!|\!|=|\!|\!|A^*B|\!|\!|$ for all $A, B\in {\mathcal A}$.
Then there exist
a unitary operator $U$ and a partial isometry $V_A$
on ${\mathcal B}(H)$,
where $V_A$ depends on $A$ and $V_A^*V_A$ is the right support projection of $A$,
such that
$$\Phi (A)=UAV_A \quad\hbox{for every}\quad A \in {\mathcal A},$$
or  
$$\Phi(A) = {U J A J V_A} \quad\hbox{for every}\quad A \in {\mathcal A}.$$
\end{corollary}

\section{$C$-numerical radius preservers}

Let $A, C\in {\mathcal B}(H)$ with $C$
being a trace class operator.
Then the $C$-numerical range and the $C$-numerical radius of $A$ are
respectively defined by$$ W_C(A)=\{\mbox{tr}(CUAU^*): U \in {\mathcal
B}(H) \ \hbox{ is unitary} \}$$
and$$w_C(A)=\sup \{|\lambda |: \lambda \in W_C(A)\}.$$
Observe that the $C$-numerical radius is unitary
similarity invariant 
and  is a semi-norm on $\cB(H)$.
In the finite dimensional case, the $C$-numerical radii can be viewed as the building blocks
of unitary similarity invariant norms $\|\cdot\|$, i.e.,
norms satisfying $\|UAU^*\| = \|A\|$ for all
$A, U \in \cB(H)$ such that $U$ is unitary, in the following sense.
For any unitary similarity invariant norm $\|\cdot\|$, there is a compact subset $\cS$ of $\cB(H)$ such that
$$\|A\| = \max\{ w_C(A): C \in \cS\}.$$
The readers may refer \cite{Li, LMR}
and the references therein for more properties of $C$-numerical
range and $C$-numerical radius.
Note also that the roles of $C$ and $A$ in the definition of $W_C(A)$
are symmetric, so $W_C(A)=W_A(C)$.

From the definition we observe also that, if $C=\sum _{i=1}^kx_i\otimes x_i$ is a rank $k$ projection,
where $\{x_i\}_{i=1}^k$ is an orthonormal set, then the $C$-numerical range and the $C$-numerical radius of $A$
reduce, respectively, to the $k$-numerical
range and the $k$-numerical radius of $A$,
$$W_k(A) =\left\{ \sum _{i=1}^k\langle Ax_i, x_i\rangle :
\{x_1, \dots, x_k \} \subseteq H \hbox{ is an orthonormal set} \right\},$$
and
$$w_k(A) = \sup \{|\lambda |: \lambda \in W_k(A) \}.$$
In particular, when $k=1$, $W_k(A)$ and $w_k(A)$ reduce to the classical 
numerical range $W(A)$ and the classical numerical radius $w(A)$ of $A$, respectively. 
One may see \cite{Halmos}
for some background of the $k$-numerical range.

If $C=x\otimes y$ with $\langle x, y\rangle =q$ and $\|x\|=\|y\|=1$,
then the $C$-numerical range and the $C$-numerical radius of $A$ reduce, respectively,
to the $q$-numerical range and the $q$-numerical radius of $A$,
$$
W_q(A)=\{\langle Ax, y\rangle :  x, y\in H\mbox{ are unit vectors such that }
\langle x, y\rangle =q\},$$
and
$$w_q(A)=\sup
\{|\lambda |: \lambda \in W_q(A)\}.$$
Clearly, $W_{qz}(A)=zW_q(A)$ for every $z\in \mathbb T$, and hence $w_{qz}(A)=w_q(A)$.
Therefore, one may assume that $q$ is real and $0 \le q \le 1$.
If $A$ is of rank one,
then $w_C(A)=w_A(C)=\|A\|w_q(C)$ with $q = \frac{|{\rm tr}(A)|}{\|A\|}$.

We observe that $C$-numerical radius satisfies conditions (F1) -- (F3) of $F(\cdot)$ stated before Theorem \ref{main}.
In fact, one can obtain the following result.

\medskip\noindent
\begin{theorem}\label{wc}
Let $H$ be a complex Hilbert space with
$\dim H\geq 3$, $\mathcal A$ and $\mathcal B$ be subsets of
${\mathcal B}(H)$ which contain all rank one operators. Assume
that $C\in {\mathcal B}(H)$ is a finite rank operator such that

\begin{enumerate}
\item[\rm (1)] $w_0(C) \ne w_1(C)$, or

\medskip
\item[\rm (2)]
$w_0(C) = w_1(C)$ and $w_q(C) > w_0(C)$ for all $q \in (0,1)$.
\end{enumerate}

\medskip
\noindent
Suppose a surjective map $\Phi :{\mathcal A}\rightarrow {\mathcal B}$ satisfies that
$w_C(\Phi (A)^*\Phi (B))=w_C(A^*B)$ for all
$A, B\in {\mathcal A}$. Then there exist unitary
operators $U$ and $V$ on $H$
and a functional $h:{\mathcal A}\rightarrow \mathbb T$ such that
$$\Phi (A)=h(A)UAV \quad \hbox{for every}\quad  A\in {\mathcal A},$$
or
$$\Phi (A)=h(A)U\overline AV \quad \hbox{for every}\quad  A\in {\mathcal A}.$$
\end{theorem}

\medskip
For any orthonormal set $\{x, y \} \subseteq H$, we have
$w_C\left(x\otimes (qx + \sqrt{1-|q|^2} y) \right)=w_q(C)$. Thus, the condition on
$w_0(C), w_1(C), w_q(C)$, etc. can be restated in terms of
$w_C(x \otimes (q x + \sqrt{1-|q|^2} y)$
for an (any) orthnormal set $\{x,y\}$. That is,
\begin{enumerate}
\item[\rm (a)] $w_C(x\otimes y) \ne w_C(x\otimes x)$, or

\medskip
\item[\rm (b)]
$w_C(x\otimes y) = w_C(x\otimes x)$ and $w_C\left(x\otimes (qx + \sqrt{1-|q|^2} y) \right) > w_C(x\otimes x)$ for all $q \in (0,1)$.
\end{enumerate}

Note also that  $w_C(X) = w_C(\overline X)$ for all $X \in \cB(H)$
if $C$ is a rank one normal operator.
On the other hand, if $\cA$ consists of only rank one operators,
then for any  trace class  operator $C \in \cB(H)$, $w_C(X) = w_C(\overline X)$ for all $X\in \cA$.
It is known that if $C$ is unitarily similar to $\overline C$, then
$\overline{W_C(X)} = W_C(\overline X)$ for any $X \in \cB(H)$ (see \cite{Li}) and hence
$w_C(X) = w_C(\overline X)$ for any $X \in \cB(H)$.
Therefore, it can be seen that the existence of the map $A\mapsto h(A) U\overline AV$
depends on $C$ and $A$.


\medskip
By checking the proof of \cite[lemma 3.1]{CH2}, one can obtain the following result.

\begin{lemma}\label{constant}
Let $H$, $\mathcal A$ and $\mathcal B$ be just as assumptions in Theorem \ref{main}.
Assume that $F^\prime :{\mathcal A}\rightarrow [0, +\infty)$ is a unitary similarity invariant function.
Let $u\in H$ be a unit vector. If there is a vector $v\in H$ which is linearly
independent of $u$ so that $F^\prime (u\otimes h)=F^\prime (v\otimes h)$ for every $h\in H,$ then,
for all unit vectors $x$, $f\in H$, $F^\prime (x\otimes f)=F^\prime (u\otimes u)$. That is,
$x \otimes f \mapsto F^\prime (x \otimes f)$  is a constant function.
\end{lemma}

We also need the following technical lemma for $q$-numerical radius.

\begin{lemma} \label{lwq}
For any $0 < q < r \le 1$ and finite rank operator $C \in \BH$,
$$w_{q}(C) \ge \min\{ w_0(C), w_r(C)\}.$$
Furthermore, when $w_0(C) \ne w_r(C)$, the inequality is always strict.
\end{lemma}

\it Proof. \rm It is known from \cite[Theorem 2.9]{LMR} that for any $q_1,q_2\in [0,1]$, 
$$d\left(W_{q_1}(C), W_{q_2}(C) \right) \le \|C\| \sqrt{ |q_1 - q_2|^2 + 2 |q_1- q_2|},$$
where $d(\cdot,\cdot)$ is the standard Hausdorff metric. It follows that the map $q\mapsto W_q(C)$ 
is a lower semi-continuous map on $[0,1]$.
Also it is shown in \cite[Lemma 5.7]{T} that for any $q_1,q_2 \in [0,1]$, if $z_i \in W_{q_i}(C)$ for $i = 1,2$, then
$$\frac{1}{2} \left( z_1 + z_2 \right) \in  W_{\frac{1}{2}(q_1+q_2)}(C).$$
Inductively, one has
$t z_1 + (1-t) z_2 \in W_{ t q_1 + (1-t) q_2}(C)$ for all $t = {k}/{2^\ell}$, where $k$ and $\ell$ are nonnegative integers and $k \le 2^\ell$.
Since the set $\{k/2^\ell: \hbox{ $k$ and $\ell$ are nonnegative integers and $k \le 2^\ell$}\}$ is dense in the set $[0,1]$, 
by the lower semi-continuity of the map $q \mapsto W_q(C)$,
one can conclude that
\begin{eqnarray}\label{convex}
t z_1 + (1-t) z_2 \in W_{ t q_1 + (1-t) q_2}(C) \quad\hbox{for all}\quad t \in [0,1].
\end{eqnarray}
It is  also known \cite[Lemma 5.5]{T} (see also \cite{Tsing}) that
the $q$-numerical range can be written as the union of circular discs, i.e.,
$$W_q(C) = \bigcup_{x\in H,\, \|x\| = 1} \left\{ z \in \IC: |z - q x^*C x| \le \sqrt{1-q^2} \sqrt{ \|Cx\|^2 - |x^*Cx|^2 } \right\}.$$
In particular,  $W_0(C)$ is a circular disc centered at the origin with radius $w_0(C)$.

Now fixed $0 < r \le 1$. Without loss of generality, we may assume that $w_r(C) \in W_r(C)$.
Otherwise, we can replace $C$ by $e^{i\theta} C$ for some $\theta \in [0,2\pi)$.
Also as $W_0(C)$ is a circular disc centered at the origin, $w_0(C) \in W_0(C)$.
For any $q \in (0,r)$, let $z =t\, w_r(C) + (1- t)\, w_0(C)$ with $t = q/ r\in (0,1)$.
Then $z\in (0,\infty)$ and it is a convex combination of $w_0(C)$
and $w_r(C)$, and hence $z \in W_{t r + (1-t) 0}(C) = W_q(C)$ by the property (\ref{convex}). Then
$$w_q(C) \ge z 
= t\, w_r(C) + (1-t)\, w_0(C) \ge \min\{w_0(C), w_r(C)\}.$$
Finally, it can be seen that the last inequality is always strict if $w_0(C) \ne w_r(C)$.
\qed

\it Proof of Theorem \ref{wc}. \rm
Suppose $\Phi: \cA \to \cB$ is a surjective map satisfying
$w_C(\Phi(A)^*\Phi(B)) = w_C(A^*B)$ for all $A,B \in \cA$.
Notice that $\Phi$ satisfies (\ref{zero}) in Theorem \ref{main}
and $w_C(\cdot)$
satisfies conditions (F1) -- (F3).
By Theorems \ref{main} and \ref{main-2} and replacing $\Phi$ by $A\mapsto \Phi(\overline A)$, if necessary,
we may assume that there exist unitary 
operator $U$ and a map $h: H \times H \to H$
with $\|h(x,f) \| = \|f\|$ 
such that
$$\Phi(x\otimes f) = Ux \otimes h(x,f) \quad\hbox{for every} \quad x,f \in H.$$
Also, condition (F2) and Lemma \ref{2.2} ensures that $\Phi$ preserves rank one operators in both directions.
Therefore, the map $f \mapsto h(x,f)$ is surjective on $H$.
Define $h_x(f) = h(x,f)$.

Since $C$ is a finite rank operator,
there exists $D$ acting on
some finite dimensional space such that $C=D\oplus 0$ according to some space
decomposition.
Then for any unit vectors  $f,g \in H$,
$w_C(f \otimes g) = w_{|\langle f,g \rangle|}(C) = w_{|\langle f, g \rangle|}(D)$.
Now for any $x, y, f, g\in H$ with $\langle x,y\rangle \ne 0$ and $\|f\| = \|g\| = 1$,
$$|\langle x,y \rangle| w_{D}(h_x(f) \otimes h_y(g)) = w_{D}(\Phi(x\otimes f)^* \Phi(y \otimes g) )
= w_{D}( (x \otimes f)^* (y \otimes g)) = |\langle x,y \rangle| w_{D}(f \otimes g).$$
Because $\|h_x(f)\| = \|h_y(g)\| = 1$, we have
\begin{multline}\label{fg}
\hspace{1cm}w_{|\langle h_x(f), h_y(g) \rangle|}(D)
= w_{h_x(f) \otimes h_y(g)}(D) 
= w_{C}(h_x(f) \otimes h_y(g)) \\
= w_{C}(f \otimes g)
= w_{f\otimes g}(D)
= w_{|\langle f,g \rangle|}(D). \hspace{.65cm}
\end{multline}
We claim that for any $x,y, f,g \in H$ with  $\langle x,y\rangle \ne 0$ and $\|f\| = \|g\| = 1$,
\begin{eqnarray}\label{claim}
\langle f,g\rangle = 0 \quad \Longleftrightarrow \quad
\langle h_x(f), h_y(g) \rangle = 0.
\end{eqnarray}
We consider three cases.

{\bf Case 1.} Suppose $w_0(D) < w_1(D)$. Then by Lemma \ref{lwq},
for any $q \in (0,1)$,
$w_q(D) >  \min\{w_0(D),w_1(D)\} = w_0(D)$.
The claim follows readily from  equation (\ref{fg}).

\medskip
{\bf Case 2.} Suppose $w_1(D) < w_0(D)$.
As the map $q\mapsto w_q(D)$ is a continuous function on $[0,1]$,
there exists some $r \in (0,1)$ such that $w_q(D) < w_0(D)$ for all $q \in [r,1]$.
By Lemma \ref{lwq}, we have
\begin{eqnarray}\label{wq1}
w_q(D) > \min\{w_0(D), w_r(D)\} = w_r(D) \quad\hbox{for all}\quad q \in [0,r).
\end{eqnarray}
On the other hand, for any $0 < r \le q_1 < q_2 \le 1$, By Lemma \ref{lwq} and the fact that $w_{q_2}(D) < w_0(D)$,
$$w_{q_1}(D) > \min\{w_0(D),  w_{q_2}(D)\} = w_{q_2}(D).$$
Therefore, the map $q\mapsto w_q(D)$ is strictly deceasing on $[r,1]$.
Together with inequality (\ref{wq1}), one can conclude that
for any $q \in [r,1]$ and $q' \in [0,1]$,
\begin{eqnarray}\label{wq}
w_{q'}(D) = w_q(D) \quad \Longrightarrow \quad q' = q.
\end{eqnarray}
Pick an integer $m\in \mathbb N$ such that $\cos(\frac{\pi}{2^m}) > r$ and let $\theta_k = \frac{\pi}{2^{m+1-k}}$ for $k = 1,\dots,m$.
We claim that for any $x,y,f,g\in H$ with $\langle x,y\rangle \ne 0$ and $\|f\| = \|g\| = 1$,
$$|\langle h_x(f),h_x(g)\rangle| = |\langle f,g\rangle|
\quad\hbox{whenever} \quad
|\langle f, g \rangle | \in [\cos \theta_k, 1]\quad \hbox{or}\quad
|\langle h_x(f), h_x(g) \rangle | \in [\cos \theta_k, 1],$$
for $k = 1,\dots, m$.

Notice that $\cos \theta_1 = \cos (\frac{\pi}{2^m} ) \in (r,1]$. By equation (\ref{wq}),
when $|\langle f,g\rangle| \in [\cos \theta_1, 1]$ or
$|\langle h_x(f), h_y(g) \rangle | \in [\cos \theta_1,1]$,
$$w_{|\langle h_x(f),h_y(g) \rangle|}(D) = w_{|\langle f,g\rangle|}(D) \quad \Longrightarrow\quad
|\langle h_x(f),h_y(g)\rangle| = |\langle f,g\rangle|.$$
Thus, the claim holds for $k =1$. Suppose now that the claim holds for some $k = \ell$.
We show that the claim also holds for $k = \ell +1$.
Suppose $|\langle f,g\rangle| = \cos  2\phi \in [\cos \theta_{\ell+1}, \cos \theta_{\ell}]$,
then $0 < \theta_\ell \le 2 \phi\le \theta_{\ell+1}$. Assume $\langle f,g\rangle = e^{it} \cos 2 \phi$ for some $t\in \IR$.
Let $z = (\sin 2\phi)^{-1} \left( g - e^{it} (\cos 2 \phi) f\right)$.
Then $z$ is a unit vector in $H$ with $\langle f,z\rangle = 0$ and
$g = e^{it} (\cos 2 \phi) f+ (\sin 2\phi) z$. Define
$u = e^{it} (\cos \phi) f+ (\sin \phi) z$.
Then direct computations show that
$|\langle f,u \rangle| =| \langle g, u \rangle| = \cos \phi \in [\cos\theta_\ell, 1]$ as $\theta_{\ell-1} = \frac{1}{2} \theta_\ell \le \phi \le \frac{1}{2} \theta_{\ell+1} = \theta_\ell$.
Then by induction assumption, we have
$$|\langle h_x(f),h_x(u)\rangle| = |\langle h_y(g),h_x(u)\rangle| = \cos \phi.$$
Then
$h_x(f) = e^{it_f} (\cos \phi) h_x(u) + (\sin \phi) z_f$ and
$h_y(g) = e^{it_g} (\cos \phi) h_x(u) + (\sin \phi) z_g$
for some $t_f,t_g\in \IR$ and unit vectors $z_f,z_g\in H$ with
$\langle h_x(u), z_f \rangle = \langle h_x(u), z_g \rangle = 0$.
It follows that
$$|\langle h_x(f),h_y(g)\rangle|
= \left| e^{i(t_g - t_f)} \cos^2 \phi  + \langle z_f, z_g\rangle \sin^2 \phi \right|
\ge \cos^2\phi - \sin^2\phi  = \cos 2\phi = |\langle f,g\rangle|.$$

Suppose $|\langle h_x(f),h_y(g)\rangle| > \cos 2\phi$. Let
$\langle h_x(f),h_y(g)\rangle = e^{is} \cos 2 \psi$ with $0 \le \psi < \phi \le \theta_\ell$.
If $\psi \le \theta_\ell$, then $\cos \psi \in [\cos \theta_\ell, 1]$, then by the claim holds by induction.
Now suppose $\theta_\ell < \psi < \theta_{\ell+1}$.
Let $d = (\sin 2\psi)^{-1} \left( h_y(g) - e^{is} (\cos 2 \psi) h_x(f)\right)$.
Then $d$ is a unit vector in $H$ with $\langle h_x(f),d\rangle = 0$ and
$h_y(g) = e^{it} (\cos 2 \psi) h_x(f)+ (\sin 2\psi) d$. Define
$w = e^{is} (\cos \psi) h_x(f)+ (\sin \psi) d$ and by the fact that $h_x$ is surjective, there exists a unit vector $v \in H$ such that
$h_x(v) =w$. Notice that $|\langle h_x(f),h_x(v)\rangle| = |\langle h_y(g),h_x(v)\rangle| = \cos \psi$.
By the induction assumption, $|\langle f,v\rangle | = |\langle g,v\rangle| = \cos \psi$.
By a similar argument as above, one can conclude that
$$|\langle f,g\rangle| \ge \cos 2\psi = |\langle h_x(f),h_y(g)\rangle|.$$
Combining the above two inequalities, we conclude that
$|\langle f,g\rangle|  = |\langle h_x(f),h_y(g)\rangle|$ if $|\langle f,g \rangle| \in [\cos \theta_{\ell+1},1]$.
The case when $|\langle h_x(f),h_y(g) \rangle| \in [\cos \theta_{\ell+1},1]$ can be proved by a similar argument.
Therefore, the claim holds and so as condition (\ref{claim}).

\medskip
{\bf Case 3.}
Assume $w_0(D) = w_1(D)$. Then by condition (2) in Theorem \ref{wc} ,  $w_q(D) > w_0(D)$ for all $q \in (0,1)$. With equation (\ref{fg}),
$$|\langle f, g \rangle | \in \{0,1\} \quad \Longleftrightarrow
\quad |\langle h_x(f), h_y(g)\rangle| \in \{0,1\}.$$
Suppose there exist
$x,y,f,g \in H$ with $\langle x,y\rangle \ne 0$ and $\|f\| = \|g\| = 1$ such that
$|\langle f,g \rangle | = 0$ and $|\langle h_x(f), h_y(g) \rangle|  = 1$.
Then $h_y(g) = \mu h_x(f)$ for some complex $\mu\in \bT$.
Pick a unit vector $\ell \in H$ such that $\langle \ell,f\rangle = \langle \ell, g\rangle = 0$.
Let $u = \alpha f + \beta g + \gamma \ell$ with
$\alpha, \beta, \gamma \ge 0$ and $\alpha^2 + \beta^2 + \gamma^2 = 1$. Then
$$w_{\alpha}(D)
= w_{|\langle f,u\rangle|}(D)
= w_{|\langle h_x(f),h_x(u) \rangle|}(D)
= w_{|\langle h_y(g),h_x(u) \rangle|}(D)
= w_{|\langle g,u \rangle|}(D)
= w_{\beta}(D).
$$
Thus, $w_\alpha(D) = w_\beta(D)$ for any $\alpha, \beta$ with $\alpha^2 + \beta^2 \le 1$.
It follows that $w_\alpha(D) = w_\beta(D)$ for any $\alpha,\beta \in [0,1]$.
But this contradicts to the condition (2) in Theorem \ref{wc}.

Now suppose there exist $x,y, f,g \in H$ with $\langle x,y\rangle \ne 0$ and $\|f\| = \|g\| = 1$
such that
$|\langle h_x(f), h_y(g) \rangle|  = 0$ and $|\langle f,g \rangle | = 1$.
Then $g = \mu f$ for some complex $\mu\in \bT$.
Pick a unit vector $\ell \in H$ such that $\langle \ell,h_x(f)\rangle = \langle \ell, h_y(g)\rangle = 0$.
Let $v = \alpha h_x(f) + \beta h_y(g) + \gamma \ell$ with
$\alpha, \beta, \gamma \ge 0$ and $\alpha^2 + \beta^2 + \gamma^2 = 1$.
Since $h_x$ is surjective, there exists a unit vector $u \in H$ such that $h_x(u) = v$.
Then
$$w_{\alpha}(D)
= w_{|\langle h_x(f),h_x(u) \rangle|}(D)
= w_{|\langle f,u\rangle|}(D)
= w_{|\langle g,u\rangle|}(D)
= w_{|\langle h_y(g),h_x(u) \rangle|}(D)
= w_{\beta}(D).
$$
This is again impossible as explained before.
Combining with the first observation,
we conclude that
 $\langle h_x(f),h_y(g)\rangle = 0$ if and only if $\langle f,g \rangle = 0$
and $|\langle h_x(f),h_y(g)\rangle| = 1$ if and only if $|\langle f,g \rangle| = 1$.
Thus, the claim (\ref{claim}) holds.

\medskip
Finally by condition (\ref{claim}) and the Uhlhorn's theorem \cite{U},
there exists a unitary or conjugate unitary operator $V_x$
and a unit complex $\mu_{x,f} \in \mathbb T$ such that
such that $h(x,f) = h_x(f) = \mu_{x,f} V_x f$ for all $f\in H$ with $\|f\| = 1$.
Suppose $V_x^*V_y$ is not a scalar for some $x,y\in H$ with $\langle x,y \rangle \ne 0$.
Then there exists a unit vector $g\in H$ such that
$g$ and $V_x^*V_y g$ are linearly independent. In this case, one can always find another unit vector $f\in H$
such that $\langle f, g\rangle = 0$ while $\langle f, V_x^*V_y g \rangle \ne 0$.
But then
$$\langle f, g\rangle = 0 \quad\Longrightarrow\quad 0 = \langle h(x,f), h(y,g) \rangle = \langle V_x f, V_y g \rangle = \langle f , V_x^*V_y g\rangle \ne 0.$$
Therefore, $V_x^* V_y$ is a scalar and hence $V_x$ and $V_y$ are linearly dependent for all $x,y\in H$
with $\langle x,y\rangle \ne 0$. Finally, when $\langle x , y \rangle  = 0$, then there always exists $z \in H$
such that both $\langle x,z \rangle$ and $\langle y,z\rangle$ are nonzero. Then $V_x$ and $V_z$ are linearly dependent
and $V_y$ and $V_z$ are also linearly dependent. Therefore, $V_x$ and $V_y$ are linearly dependent.
In this case, we can conclude that $h(x,f) = \mu_{x,f} V f$ for some $V$.
Then there exists a functional $d:{\mathcal F}_1(H)\rightarrow
{\mathbb T}$ such that $\Phi (x\otimes f)=d(x\otimes f)Ux\otimes Vf$ for
all $x, f\in H$. Since $U$ is linear, it follows that $V$ must be linear (and not conjugate linear) too.
Now for
every $A\in {\mathcal A}$ with the rank greater than one and all $x, f\in H$, we have
\begin{eqnarray}\label{h}
w_C( (V \Phi (A)^*
U) x\otimes f)
=w_C(\Phi (A)^*\Psi (x\otimes f))
=w_C(A^*x\otimes
f).
\end{eqnarray}
Then Lemma \ref{constant} implies that $V\Phi(A)^*U x$ and $A^*x$ are linearly dependent, it follows from \cite{Kap} (see also \cite{Hou}) that locally linearly dependent implies linearly dependent,we have
$V\Phi(A)^*U$ and $A^*$ are linearly dependent, and hence
$U^*\Phi (A)V^*$ and $A$ are linearly dependent. So, there exists a functional $h $ on ${\mathcal
A}$ such that $\Phi (A)=h (A)UAV$ for every $A\in {\mathcal A}$, where $h (A)=d(A)$ if $A$ is of rank one.
Equation (\ref{h}) ensures that $|h (A)|=1$ for every $A\in {\mathcal
A}$. The result follows.
\qed


\medskip\noindent
{\bf Remark 4.4}. It would be nice to prove Theorem \ref{wc} without the assumptions (1)--(2) on $C$.

\section{pseudo spectral radius preservers}

For every $\varepsilon >0$, define  the $\varepsilon$-pseudo spectrum
$\sigma_\varepsilon (A)$ and the
 $\varepsilon$-pseudo spectral radius $r_\varepsilon (A)$ as
$$
\sigma _\varepsilon (A)=\{z\in {\mathbb C}:
\|(zI-A)^{-1}\|>\varepsilon ^{-1}\},$$
and
$$r_\varepsilon (A)=\sup \{|z|:z\in \sigma_\varepsilon (A)\},$$
respectively.

From the definition, it follows that the $\varepsilon$-pseudo spectra of $A$
are a family of strictly nested closed sets, which grow to fill
the whole complex plane as $\varepsilon \rightarrow \infty$, and
that the intersection of all the pseudo spectra is the spectrum,
$$\bigcap _{\varepsilon >0}\sigma _\varepsilon (A)=\sigma (A),$$
where $\sigma (A)$ denotes the spectrum of
$A$.  Accordingly, we have
$$\lim_{\varepsilon \rightarrow 0^+}r_\varepsilon (A)=r(A),$$
where $r(A)$ is the spectral radius of $A$.

In the following,
we will apply Theorems \ref{main} and \ref{main-2} to
characterize maps preserving the pseudo spectral radius.
Remark that the pseudo spectrum preservers of standard product 
$AB$ and Jordan product $AB+BA$ were determined in \cite{CL1}.

\begin{theorem}\label{psr}
Let $H$, $\mathcal A$ and $\mathcal B$ be
just as assumptions in Theorem \ref{main} and $\varepsilon >0$. 
Suppose a surjective map $\Phi:{\mathcal A}\rightarrow \mathcal B$ 
satisfies that 
$$r_\varepsilon (\Phi (A)^*\Phi (B))=r _\varepsilon
(A^*B)\quad \mbox{ for all}\quad A, B\in \mathcal A.$$ 
Then there exist unitary operators $U$ and $V$ on $H$
and a functional $h:{\mathcal A}\rightarrow \mathbb T$ such that
$$\Phi (A)=h(A)UAV \quad\hbox{for every}\quad A\in {\mathcal A},$$
or 
$$
\Phi (A)=h(A)U\overline AV
\quad\hbox{for every}\quad A\in {\mathcal A}.$$
\end{theorem}

We need some preliminary results to prove Theorem \ref{psr}.
First, it is clear from definition that the pseudo spectrum and
pseudo spectral radius are unitary
similarity invariant. Now let us recall other properties of the
pseudo spectrum (see \cite{TE1, TE}). Let $\varepsilon > 0$ be arbitrary
and $D(a, \varepsilon )= \{ \mu \in {\mathbb C} : |\mu - a| <
\varepsilon \}$, where $a \in \mathbb C$. For $A\in {\mathcal B}(H)$,
\begin{enumerate}
\item[\rm (i)]  $\sigma (A) + D(0,\varepsilon ) \subseteq \sigma _\varepsilon
(A)$.

\item[\rm (ii)] If $A$ is normal, then $\sigma _\varepsilon (A)=\sigma (A) +
D(0,\varepsilon ) $.

\item[\rm (iii)] For any $c\in \mathbb C$, $\sigma _\varepsilon
(A+cI)=c+\sigma _\varepsilon (A)$.

\item[\rm (iv)] For any nonzero $c\in \mathbb C$, $\sigma _{\varepsilon
}(cA)=c\,\sigma _{\frac{\varepsilon }{|c|}}(A)$.
\end{enumerate}
The following lemmas were proved in \cite{CL1}.
In the finite dimensional case, one sees also \cite{CL}.

\begin{lemma}\label{4.1}
Let $A\in
{\mathcal B}(H)$, $\varepsilon >0$ and $a\in \mathbb C$. Then the followings hold.
\begin{enumerate}
\item $A=aI$ if and only if $\sigma _\varepsilon
(A)=D(a, \varepsilon )$.

\item
$A=aP$ for some nontrivial projection $P\in {\mathcal
B}(H)$ if and only if $\sigma _\varepsilon (A) = D(0,\varepsilon )
\cup D(a,\varepsilon )$.
\end{enumerate}
\end{lemma}

\begin{lemma}\label{4.2}
Let $\varepsilon >0$ and $x, f \in H$.
\begin{enumerate}
\item Then
$$r_\varepsilon (x\otimes f)=\frac{1}{2}
( \sqrt{|\langle x, f\rangle |^2 + 4\varepsilon ^2+4 \varepsilon
\|x\|\|f\|} + |\langle x, f\rangle |),$$
which is attained at a point in $\sigma_{\varepsilon}(x\otimes f)$
with the direction $\langle x, f\rangle$.

\item
$\langle x, f\rangle =0$ if and only if
$$\sigma _\varepsilon (x\otimes f)=
D(0, \sqrt{\varepsilon^2+\|x\|\|f\|\varepsilon}).$$
\end{enumerate}
\end{lemma}

Evidently, the pseudo spectral radius satisfies
conditions (F1) -- (F3) stated before Theorem \ref{main}.
We can now present the proof of Theorem \ref{psr}.

\it Proof of Theorem \ref{psr}. \rm
Notice that $\Phi$ satisfies (\ref{zero}) in Theorem \ref{main} 
and $r_\varepsilon(\cdot)$ satisfies
conditions (F1) -- (F3).
By Theorems \ref{main} and \ref{main-2}, and replacing $\Phi$ with $A\mapsto \Phi(\overline A)$, if necessary,
we may assume that
there exist a unitary operator $U$ and a map $h: H \times H \to H$
with $\|h(x,f) \| = \|f\|$ such that
$$\Phi(x\otimes f) = Ux \otimes h(x,f) \quad\hbox{for every} \quad x,f \in H.$$
Define $h_x(f) = h(x,f)$. First, for any $x,y, f,g\in H$
with $\langle x,y\rangle \ne 0$,
$$|\langle x,y \rangle| r_{\frac{\varepsilon}{|\langle x,y\rangle|}}(h_x(f) \otimes h_y(g))
= r_\varepsilon(\Phi(x\otimes f)^* \Phi(x \otimes g) )
= r_\varepsilon( (x \otimes f)^* (x \otimes g))
= |\langle x,y \rangle| r_{\frac{\varepsilon}{|\langle x,y\rangle|}}(f \otimes g),$$
and hence $r_{\hat \varepsilon}(h_x(f) \otimes h_x(g))
= r_{\hat \varepsilon}(f \otimes g)$
with $\hat \varepsilon ={\frac{\varepsilon}{|\langle x,y\rangle|}}$.
By Lemma \ref{4.2},
\begin{multline*}
\frac{1}{2} (
\sqrt{|\langle h_x(f), h_y(g)\rangle |^2 + 4\hat\varepsilon ^2+4 \hat \varepsilon
\|h_x(f)\|\|h_y(g)\|} + |\langle h_x(f), h_y(g)\rangle | ) \cr
= \frac{1}{2}(
\sqrt{|\langle f, g\rangle |^2 + 4\hat \varepsilon ^2+4 \hat \varepsilon
\|f\|\|g\|} + |\langle f, g\rangle |).
\end{multline*}
With the fact that $\|f\| = \|h_x(f)\|$ and $\|g\| = \|h_y(g)\|$, we conclude that
\begin{eqnarray}\label{claim2}
\langle f,g\rangle = 0 \quad \Longleftrightarrow \quad
\langle h_x(f), h_y(g) \rangle = 0.
\end{eqnarray}
Then using the same argument as in the last part of Theorem \ref{wc}, the result follows.
\qed

\begin{corollary}
Let  $\varepsilon >0$.
Suppose a surjective map $\Phi
:\cB(H) \rightarrow \cB(H)$ satisfies that $$\sigma
_\varepsilon (\Phi (A)^*\Phi (B))=\sigma _\varepsilon
(A^*B)\quad \mbox{ for all}\quad A, B\in \mathcal A.$$ 
Then there exist a complex unit $\alpha $ and unitary operators $U, V\in {\mathcal B}(H)$ such that
$$\Phi (A)=\alpha UAV \quad \hbox{for every}\quad A\in \cB(H).$$
\end{corollary}

\medskip
{\it Proof.} Clearly, $r_\varepsilon (\Phi (A)^*\Phi (B))=r_\varepsilon (A^*B)$ for all $A$, $B\in \mathcal A$.
Then Theorem \ref{psr} implies that there exist unitary operators $U$ and $V$ on $H$
and a functional $h:\cB(H)\rightarrow \mathbb T$ such that
$\Phi (A)=h(A)UAV$ or $\Phi(A) = h(A) U\overline A V$ for every $A\in \cB(H)$.
Suppose $\Phi(A) = h(A) U\overline A V$ for all $A \in \cB(H)$.
It follows that for
every $A\in \cB(H)$,
\begin{eqnarray}\label{bar}
\sigma _\varepsilon (A)=\sigma
_\varepsilon (\Phi(I)^* \Phi (A))=\sigma _\varepsilon (\overline{h(I)}h(A)V^*\overline AV)
= \sigma _\varepsilon (\overline{h(I)}h(A)\overline A).
\end{eqnarray}
Take an arbitrary nontrivial projection $P\in \cB(H)$ and let $A=i P+(1-i) (I-P)$.
Then Equation (\ref{bar}) and Lemma \ref{4.1} imply that
$$D(\overline{h(I)}h(A)(-i), \varepsilon )\cup
D(\overline{h(I)}h(A)(1+i), \varepsilon
)=D(i , \varepsilon )\cup D(1-i , \varepsilon ),$$it follows
that either $\overline{h(I)}h(A)(-i)=i$ and
$\overline{h(I)}h(A)(1+i)=1-i$ or
$\overline{h(I)}h(A)(-i)=1-i$ and
$\overline{h(I)}h(A)(1+i)=i$, which is a contradiction.
Therefore, $\Phi$ has the form $A \mapsto  h(A) UAV$. Let $\Psi (A)=\Phi (I)^*\Phi
(A)$ for every $A\in \cB(H)$. Since $\Phi (I)=h(I)UV$
is unitary, we have $\Psi (I)=I$ and $\sigma _\varepsilon (\Psi
(A)^*\Psi (B))=\sigma _\varepsilon (A^*B)$ for all $A, B\in
\cB(H)$, also $\sigma _\varepsilon (\Psi (A))=\sigma
_\varepsilon (A)$ for every $A\in \cB(H)$.
We claim that
$h(A)$ is a complex unit which is independent of $A$.
Now a similar discussion just as in Equation (\ref{bar}) implies that $$\sigma _\varepsilon (A)=\sigma _\varepsilon
(\overline{h(I)}h(A)A)=\overline{h(I)}h(A)\sigma_\varepsilon
(A)\quad \mbox{ for all }A\in \cB(H),$$
which implies that, for every $A\in \cB(H)$,
$$z\in \sigma _\varepsilon (A)\Leftrightarrow \overline{h(I)}h(A)z\in \sigma _\varepsilon (A).$$
If  $\sigma_\varepsilon(A)$ is not a circular disk with center zero, then
$\overline{h(I)}h(A)=1$, and hence $h(A)=h(I)$. Now assume that $\sigma_\varepsilon(A)$
is a circular disk with center zero. Take any $x, f\in H$ satisfying $\langle x, f\rangle \neq 0$ and
$\langle A^*x, f\rangle \neq 0$. Then $h(x\otimes f)=h(I)$ and $\sigma _\varepsilon (A^*x\otimes f)$ is
not a circular disk with center zero. Thus $\sigma _\varepsilon (A^*x\otimes f)
=\sigma _\varepsilon (\Phi (A)^*\Phi (x\otimes f))
=\overline{h(A)}h(x\otimes f)\sigma _\varepsilon (A^*x\otimes f)$ implies that
$\overline{h(A)}h(x\otimes f)=1$, and hence $h(A)=h(x\otimes f) = h(I)$.
Therefore for every $A\in \mathcal A$, we have $h(A)=h(I)$.
This completes the proof.
\qed

\bigskip\noindent
{\bf Acknowledgment}

The authors would like to thank the referee for some helpful suggestions
and drawing reference \cite{Kap} to their attention.

\end{document}